\documentclass{article}
\usepackage{amsmath}
\usepackage{amscd}
\usepackage{latexsym}
\usepackage{amssymb}

\title{Graphs, Frobenius functionals, and the classical Yang-Baxter
equation}
\author{M. GERSTENHABER and A. GIAQUINTO}

\begin{document}
\maketitle
\newtheorem{Theorem}{Theorem}
\newtheorem{Corollary}{Corollary}
\newtheorem{Lemma}{Lemma}
\renewcommand{\abstractname}{}
\newcommand{\mn}{\ensuremath{\operatorname{mn}}}
\newcommand{\cn}{\ensuremath{\operatorname{cn}}}
\newcommand{\spl}{\ensuremath{\operatorname{sl}}}
\newcommand{\Spl}{\ensuremath{\operatorname{Sl}}}
\newcommand{\rk}{\ensuremath{\operatorname{rk}}}
\newcommand{\ad}{\ensuremath{\operatorname{ad}}}
\newcommand{\tr}{\ensuremath{\operatorname{tr}}}
\newcommand{\vx}{\ensuremath{\operatorname{vx}}}
\newcommand{\idx}{\ensuremath{\operatorname{idx}}}
\newcommand{\ot}{\ensuremath{\leftarrow}}
\newcommand {\f}{{\frak f}}
\newcommand {\g}{{\frak g}}
\newcommand {\n}{{\frak n}}
\newcommand {\h}{{\frak h}}
\newcommand {\PC}{{\mathcal{P}}}
\newcommand{\G}{{\Gamma}}
\newcommand{\e}{{\varepsilon}}
\newcommand{\fo}{\frak{f}_{\textrm{odd}}}
\newcommand{\fe}{\frak{f}_{\textrm{even}}}
\newcommand{\bsq}{{\blacksquare}}
\newcommand{\sq}{{\square}}
\newcommand{\bt}{{\bullet}}
\newcommand{\btl}{{\blacktriangle}}
\newcommand{\tl}{{\triangle}}
\newcommand{\btr}{{\blacktriangleright}}
\newcommand{\tgr}{{\triangleright}}
\vspace{-7mm}
\date{}
{\noindent ${}^1$\textit{Department of Mathematics, University of
Pennsylvania, Philadelphia, PA 19104-6395, \linebreak[0]
U.S.A.\linebreak[0]
email:mgersten@math.upenn.edu}\\
\noindent ${}^2$\textit{Department of Mathematics, Loyola University
Chicago, Chicago, IL 60626-5385 U.S.A.,
email:tonyg@math.luc.edu}}

\begin{abstract}\noindent \textbf{Abstract.} A Lie algebra
$\g$ is Frobenius if it admits a linear functional $F\in \g^*$ such
that the Kirillov form $B_F(x,y) = F([x,y])$ is non-degenerate. If
$\g$ is the $m$th maximal parabolic subalgebra $\PC(n,m)$ of
$\spl(n)$ this occurs precisely when $(n,m) = 1$. We define a
\emph{cyclic} functional $F$ on $\PC(n,m)$ and prove it is
non-degenerate using properties of graphs associated to $F$. These
graphs also provide in certain cases readily computable associated
solutions of the classical Yang-Baxter equation. We define the
\emph{full local ring} of a graph from which we show that the graph
can be reconstructed (as well as a \emph{reduced} local ring),
consider the seaweed algebras of Dergachev and Kirillov, and examine
the degeneration of solutions to the modified classical Yang-Baxter
equation.
\end{abstract}

\section{Introduction}

Let $\g$ be a finite dimensional Lie algebra over a field $K$ and
$F\in \g^*$ be a linear functional. The associated skew bilinear
\emph{Kirillov form} $B_F$ is defined by $B_F(x,y) = F([x,y])$ for
$x,y\in\g$. In this paper $K$ will have characteristic zero although
some important assertions hold more generally. The \emph{index} of
$\g$, $\idx(\g)$, is the minimum dimension of $\ker B_F$ as $F$
ranges over $\g^*$. Those $F$ for which the minimum is achieved are
called \emph{regular} and form a Zariski open and dense subset
$\g^*_{\mathrm{reg}}$ of $\g^*$. Clearly $\g$ operates on
$\g^*_{\mathrm{reg}}$; if $F\in \g^*, x\in \g$ then $[x,F](y) =
F([x,y])$ for all $y \in \g$.   The Lie algebra $\g$ is
\emph{Frobenius} if its index is zero, i.e., if there exists an
$F\in\g^*$ such that $B_F$ is non-degenerate. We will call such an
$F$ a \emph{Frobenius functional} and denote Frobenius Lie algebras
by $\f$.

Frobenius Lie algebras are intimately connected to skew solutions of
the classical Yang-Baxter equation (CYBE). An element $r\in \g\wedge
\g$ is a solution to the CYBE if $[r,r]=0$, where $[-,-]$ is the
Schouten bracket on $\bigwedge ^*\g$. For any Lie algebra $\g$ and
$F\in\g^*$  the bilinear form $B_F$ is by definition a coboundary in
the Chevalley-Eilenberg complex. Belavin and Drinfel'd call $\g$
\emph{quasi-Frobenius} if it admits a non-degenerate skew bilinear
form $B$ (not necessarily of the form $B_F$) which is a $2-$cocycle.
For such Lie algebras they note that if $B_{ij}$ is the matrix of
$B$ relative to some basis $x_1,\dots,x_d$ of $\f$ then
$\sum(B^{-1})_{ij}x_i\wedge x_j$ is a skew solution to the classical
Yang-Baxter equation (CYBE), \cite{BelDrin:2}, \cite{BelDrin:1}.

For applications, it is useful to have explicit functionals $F$ so
that the $r$-matrix can be exhibited. Our main interest here is
$\PC(n,m)$, the $m$th maximal parabolic subalgebra of $\spl(n)$ (the
set of all $n \times n$ matrices $X$ of trace zero with $x_{ij} = 0$
whenever both $i \ge m+1$ and $j \le m$).  A theorem of Elashvili
asserts in particular that it is Frobenius if and only if $(m,n) =
1$, \cite{Elash:1982} (where it is noted that a partial result
appears in Rais \cite{Rais:1978}, but this part covers the case here
only of $m=1$.) While any generic $F \in\PC(n,m)^*$ will be
Frobenius, such functionals are not feasible for computations.

In this note we use graphical methods to show that a canonically
defined ``cyclic'' functional on $\PC(n,m)$ is Frobenius, and that
in many cases the associated $r$-matrix is readily computable; these
results are proved in Sections \ref{sec:gamma}--\ref{sec:CYBEgraph}.
The cyclic functional is determined by writing the integers
$\{1,2,\dots,n\}$ in the cyclic order $\{1,m+1,2m+1,\dots
,(n-1)m+1\}$, where each element greater than $m$ is reduced modulo
$m$. Our computations make substantial use of the ``principal
element'' of $\f$, which was introduced in \cite{G:principal} and is
reviewed in Section \ref{sec:principal}. Other examples of Frobenius
functionals on $\PC(n,m)$ are discussed, as well as the ``seaweed''
algebras of Dergachev-Kirillov \cite{DergKir:Index} (called
``biparabolic'' by A. Joseph \cite{Joseph:biparaI}); see Sections
\ref{sec:functionals}--\ref{sec:seaweed}.

In Section \ref{sec:MCYBE} we revisit some comments of
\cite{GGS:Construction} pertaining to degenerations of solutions to
the modified classical Yang-Baxter equation (MCYBE). An element
$r\in \g\wedge \g$ is a solution to the MCYBE if $[r,r]$ is a
non-zero invariant element of $\g\wedge \g \wedge \g$. In
\cite{GGS:Construction} we stated without proof a remark about the
moduli space of solutions to the MCYBE which asserted, in effect,
which ones were limits of others in that they lie in boundaries of
their orbits, but the proof was omitted. This is quite easy to
demonstrate for $\spl(n)$ and we do so here. There seems, however,
to be a relation between this and properties of the principal
element associated to the ``cyclic'' functional which we define.
Lastly, in Section \ref{sec:local}, we define two local rings
associated to a graph which, although not used in our main results,
might be of independent interest. The ``full'' local ring enjoys the
property that it completely characterizes the graph, with just one
exception.

\section{The principal element}\label{sec:principal}

In this section we recall some results of \cite{G:principal} which
will be needed throughout the paper. Let $F$ be a Frobenius
functional on a Lie algebra $\f$. The natural map $\f\to\f^*$
defined by $x \mapsto F([x,-])$ is then invertible; the image of $F$
under the inverse is called the \emph{principal element} of $\f$ and
will be denoted $\hat F$. It is the unique element of $\f$ such that
$F([\hat F,x]) = F(x)$, or $F\circ \ad \hat F = F$; it depends, of
course, on the choice of Frobenius functional. Let $\mathfrak G$ be
the adjoint algebraic group of $\f$, i.e., the smallest algebraic
Lie group whose Lie algebra contains $ \ad\g$. Then $\frak G$
operates on $\f^*$ and the subset of Frobenius functionals is stable
under this action but we will see that the action need not be
transitive.

When $\f$ is a subalgebra of a simple Lie algebra $\g$ we will say
that $\f$ is \emph{saturated} if it is not an ideal of any larger
subalgebra of $\g$ (e.g., if it contains a Cartan subalgebra of
$\g$). In that case, $\hat F$ is semisimple, see Theorem 1 of
\cite{G:principal}. (It need not be in general; for a simple example
cf. Ooms \cite{Ooms:1980}.) For $\g = \spl(n)$, Theorem 2 of
\cite{G:principal} asserts that the eigenvalues of $\ad \hat F$ must
all be integers and independent of the choice of Frobenius
functional $F$. (Dergachev has communicated to the authors that this
holds more generally for algebraic $\g$.) The proof of Theorem 2 in
\cite{G:principal} actually shows that in this case the eigenvalues
of $\hat F$ are integers and constitute a single unbroken string,
i.e., if $i < j$ are eigenvalues then so is any $k$ with $i <k <j$.
Under the conditions above (which are probably too restrictive) one
sees that no eigenvalue of $\hat F$ can be larger than $n$, so
the eigenvalues are bounded and there are, up to similarity, only a
finite number of possibilities for $\hat F$. In a communication to
the authors Dergachev has shown that the eigenvalues of $\hat F$ do
not determine $\f$ (contrary to a conjecture in \cite{G:principal})
but it may still be the case that for any given $n$ there are only a
finite number of Frobenius subalgebras of $\spl(n)$. Denoting the
$\lambda$-eigenspace of $\f$ by $f_{\lambda}$ note that since $\dim
\f_{\lambda} = \dim \f_{1-\lambda}$ one has also that $\tr(\ad \hat
F) = \frac{1}{2}\dim \f$, \cite{Ooms:1980}. A method of calculating
the principal element for certain Frobenius fuctionals will be given
in the next section.

\section{The graph $\gamma(S)$}\label{sec:gamma}When
$\f$ is Frobenius the eigenvalues of $\hat F$ can be computed from
any regular $F \in\f^*$ (in particular, from the generic element),
but for $\f \subset \spl(n)$ it is most convenient to do so from a
``small'' Frobenius functional \cite{G:principal}. Let $e_{ij}$
denote the $n \times n$ matrix with $1$ in the $(i,j)$ place and
zeros elsewhere. If $S$ is a subset of the indices $(i,j), i\le i,j
\le n$ then $F_S$ will denote the functional $\sum_{s\in S}e_s^*$.
It is defined on the space $M_n$ of all $n \times n$ matrices but
will tacitly be restricted, without change in notation, to any Lie
subalgebra $\g$ of $M_n$ we are considering, it being understood
then that those $e_{ij}$ with $(i,j) \in S$ lie in $\g$. We will say
that $S$ \emph{carries} or \emph{supports} $F_S$. The \emph{directed
graph of the functional}, $\gamma(S)$, has vertices the integers
$1,\dots,n$ with an arrow from $i$ to $j$ whenever $(i,j) \in S$. We
will call $F_S$ \emph{small} if $\gamma(S)$ is a tree (by definition
connected, hence having all the integers \,$1,\dots, n$\, as
vertices and having exactly $n-1$ links). One can show, again by the
arguments of \cite{G:principal}, that for the maximal parabolic
subalgebras of $\spl(n)$ if $\#S < n-1$ then $F_S$ can not be
Frobenius.

An important class of functionals on the seaweed subalgebras $\g$ of
$\spl(n)$ was given by Dergachev and Kirillov \cite{DergKir:Index}.
These subalgebras, which include the maximal parabolic ones, are
discussed later. The Dergachev-Kirillov functionals have the form
$F_S$ for some $S$ and are always regular; when $\g$ is Frobenius
they are small Frobenius functionals. The `meander' introduced in
\cite{DergKir:Index} is just $\gamma(S)$. For the maximal parabolic
subalgebras of $\spl(n)$ we construct some other Frobenius
functionals, principally the ``cyclic'' functionals, for which the
associated solution of the CYBE can sometimes be effectively
calculated. This will involve some elementary properties of graphs,
next section. We conjecture that every saturated Frobenius
subalgebra of $\spl(n)$ has a small Frobenius functional.

Suppose again that we have a small linear functional $F_S$ on
$\spl(n)$ supported by a set $S$.  As in \cite{G:principal}, set
$\e_i = e_{ii} - 1/n$. These all have trace zero, one has
$\e_1+\e_2+\cdots+\e_n = 0$, and any $n-1$ of them will serve as a
basis for the Cartan subalgebra $\h$; we generally use
$\e_1,\dots,\e_{n-1}$. Denote by $K(S)$ the space spanned over the
field $K$ by the $e_s$ with $s\in S$.  The space spanned by the
$\e_i$ is then dual to $K(S)$ under the bilinear form $B_F$. One can
exhibit explicitly the dual basis to the $e_s, s \in S$: Note that
if $s = (i,j) \in S$ then removing the edge $i\to j$ disconnects
$\gamma(S)$ so every vertex remains connected precisely to one of
$i$ and $j$. Let $d_s$ be the sum of all those $\e_k$ where $k$
remains connected to $i$ or equivalently, the negative of the sum
where $k$ remains connected to $j$. (If $n$ is connected to $i$, the
former will involve $\e_n$ and the letter will not; similarly if $n$
is connected to $j$.) Then $B(d_s, e_{s'}) = \delta_{s,s'}$,\,
\cite{G:principal}. The $d_s$ are linearly independent but somewhat
more is the case. The directed graph $\gamma_S$ defines a partial
order on the set $\{1,\dots,n\}$. Conjugating by a suitable
permutation matrix we may assume that $n$ is a terminal vertex of
$\gamma_S$ and that the partial order is compatible with the natural
order. The $d_s$, which can now simply be numbered as
$d_1,\dots,d_{n-1}$, by their construction have the property that
each $d_i$ is a linear combination only of those $\e_j$ with $j\le
i$, with the coefficient of $\e_i$ equal to 1. Thus the linear
transformation giving the $d_s$ in terms of $\e_1, \dots, \e_{n-1}$
in fact has determinant equal to 1. The $d_s$ span the Cartan
subalgebra $\h$ of diagonal traceless matrices of $\spl(n)$. Set
$D_S = \sum_{s\in S}d_s$. Then one has $[D_S, e_s] = e_s$ for all $s
\in S$. The eigenvalues of $\ad(D_S)$ on $M_n$ (hence also on
$\spl(n)$) are necessarily integers.

To illustrate this, below in diagrammatic form are two small
Frobenius functionals on $\PC(7,3)$, the ``cyclic'' one of \S
\,\ref{sec:cyclic} (where it will be shown to be Frobenius), and
that of Dergachev-Kirillov (which is known to be Frobenius,
\cite{DergKir:Index}). In the two diagrams ``$\mathbf{x}$'' marks
the matrix entries which support the functional, light dots mark the
places where the matrix entries must be zero, and the diagonal is
marked visually by dark dots. The Dergachev-Kirillov functional is
constructed by placing $\mathbf{x}$\,s on antidiagonals (lines along
which $i+j$ is constant) starting at the corners and proceeding
until one reaches the main diagonal. The rank of
$S_{\mathrm{cyclic}}$ (replace each $\mathbf{x}$ by $1$ and all
other entries by $0$) is four while that of $S_{\mathrm{DK}}$ is
five, so they cannot be conjugates.

\begin{equation*}
S_{\mathrm{cyclic}} = \begin{pmatrix}
\bullet&\circ&\circ&\mathbf{x}& \circ&\circ&\circ\\
\circ&\bullet&\mathbf{x}&\circ&\mathbf{x}&\circ&\circ\\
\mathbf{x}&\circ&\bullet&\circ&\circ&\mathbf{x}&\circ\\
\cdot&\cdot&\cdot&\bullet&\circ&\circ&\mathbf{x} \\
\cdot&\cdot&\cdot&\circ&\bullet&\circ&\circ \\
\cdot&\cdot&\cdot&\circ&\circ&\bullet&\circ \\
\cdot&\cdot&\cdot&\circ&\circ&\circ&\bullet
\end{pmatrix},\quad
S_{\mathrm{DK}} =\begin{pmatrix}
\bullet&\circ&\circ& \circ&\circ&\circ&\mathbf{x}\\
\circ&\bullet&\circ&\circ&\circ&\mathbf{x}&\circ\\
\mathbf{x}&\circ&\bullet&\circ&\mathbf{x}&\circ&\circ\\
\cdot&\cdot&\cdot&\bullet&\circ&\circ&\circ \\
\cdot&\cdot&\cdot&\circ&\bullet&\circ&\circ \\
\cdot&\cdot&\cdot&\circ&\mathbf{x}&\bullet&\circ \\
\cdot&\cdot&\cdot&\mathbf{x}&\circ&\circ&\bullet
\end{pmatrix}
\end{equation*}
The graph $\gamma(S_{\mathrm{DK}})$ is a chain if one disregards the
direction of the arrows; this holds for all the Dergachev-Kirillov
functionals of \cite{DergKir:Index}.
$$\gamma(S_{\mathrm{DK}}): \quad
    2\to6\to5\ot3\to1\to7\to4$$
By contrast, the graph $\gamma(S_{\mathrm{cyclic}})$ is a rooted
tree with all arrows directed away from the root (which in this
example is $2$); this holds for all the cyclic functionals defined
in Section~\ref{sec:cyclic}.
$$ \gamma(S_{\mathrm{cyclic}}):\,
\begin{array}{ccccccc}
2&\to&3&\to&6&&\\
   \downarrow&&\downarrow&&&&\\
   5&&1&\to&4&\to&7
\end{array}$$
One can compute $D_S$ (which turns out to be the principal element
$\hat F$) and its eigenvalues directly from these graphs. Suppose we
want the coefficient of $\e_i$. Removing a single link from the
graph disconnects it and leaves $i$ either connected to or
disconnected from $n$ (which here is $7$). Removing one link at a
time, the coefficient of $\e_i$ =(\#times $i$ remains connected to
$n$)$-$(\#times $i$ is disconnected from $n$). In the examples
(calculating without the use of $\e_7$) we have
\begin{align}
 D_{S_{\mathrm{DK}}}& = \operatorname{diag}(1,\ 3,\ 2,-1,\ 1,\ 2,\ 0)
-\ (8/7)\,\mathbf{I}_7\\ D_{S_{\mathrm{cyclic}}}&=
\operatorname{diag}(\ 2,\ 4,\ 3,\ 1,\ 3,\ 2,\
0)-(15/7)\,\mathbf{I}_7\label{eqn:D_S}
\end{align}
where $\mathbf{I}_7$ is the $7 \times 7$ unit matrix. These are
conjugate within the parabolic subgroup of the special linear group
$\operatorname{SL}(7)$ corresponding to removal of the $3$rd
negative root. All the $e_{ij}, i,j \le n$ are eigenvectors for $\ad
D_S$. To calculate the eigenvalue from $\gamma(S)$, note that
there is a unique path on the graph from $i$ to $j$; the eigenvalue
is (\# arrows traversed in the direction of the arrow) $-$ (\#
arrows traversed in the reverse direction). The eigenvalues (counted
with their multiplicities) of $D_{S_{\mathrm{DK}}}$ and
$D_{S_{\mathrm{cyclic}}}$ are obviously the same; more important,
they are the same on the subset consisting of those $e_{ij}$ which
lie in $\PC(7,3)$, as we know must be the case in general. (An
alternative and frequently easier calculation of $D_S$: build a
diagonal matrix $\operatorname{diag}(c_1,\dots,c_n)$ by setting $c_1
= 0$ and defining the remaining entries by requiring that $c_i-c_j =
1$ whenever $i\to j$ in $\gamma(S)$, then subtract a suitable
multiple of the identity to reduce the trace to zero.)

When $F$ is Frobenius $D_S$ is its principal element $\hat F$,
\cite{G:principal}. One can verify here directly Ooms' observation
that the sum of its eigenvalues on $\PC(7,3)$ must be $(1/2)\dim
\PC(7,3) =18$. Since the eigenvalue on $e_{ij}$ is the negative of
that on $e_{ji}$ it is sufficient to sum the eigenvalues on the 12
of those $e_{ij} \in \PC(7,3)$ whose transposes are not in
$\PC(7,3)$. Using the cyclic functional and going by rows, these are
$1, -1, 0, 2; 3, 1, 2, 4; 2, 0, 1, 3$.

Although we assume throughout that $K$ is a field of characteristic
zero, in some places
it is sufficient that it be a commutative unital ring where, as in the discussion
above, when considering $\spl(n)$ one must assume further that $1/n$
is invertible.

\section{Matching number and index of a graph}To prove that the
functionals we define are Frobenius it is useful to have some
elementary observations about graphs. (There are many treatises; a
brief review of the concepts relevant here, with references, can be
found, for example, in the Wikipedia article \cite{Wiki:match} and
\cite{Pulleyblank}.) All our graphs will be finite, not necessarily
connected, with any two vertices (or nodes) joined by at most one
edge (or link), and no edge going from a vertex to itself. Edges
will be called disjoint if they do not share a common vertex. A
\emph{matching} in a graph $\G$ is a set of disjoint edges. A
maximal matching is one which can not be enlarged, but these need
not all have the same number of edges. A \emph{maximum} matching is
one having the largest possible number of edges; this number is
called the \emph{matching number} of the graph and will be denoted
$\mn(\G)$. A \emph{perfect} or \emph{complete} matching is one which
covers every vertex of the graph, i.e., such that every vertex is an
end of some edge in the matching. Perfect matchings are necessarily
maximum. Computing matching numbers is a basic problem in graph
theory but in the case of a tree or forest (a disjoint union of
trees) there is a simple algorithm (a trivial case of some more
sophisticated ones, cf. e.g. \cite{Pulleyblank} but all we need).
Call a vertex at which more than two edges meet a branch point and
one met by only one edge an end; a terminal vertex will be one which
is either an end or is isolated (not met by any edge). A graph will
be called a \emph{chain} if it is connected, has no branch points,
and is not an isolated point. Its length is the number of edges.
When graphs are directed we can distinguish beginning or ``initial"
and end or ``terminal" vertices but an initial vertex will still be
considered terminal in the sense of undirected graphs. A
\emph{terminal chain} is a subgraph which is a chain one end of
which is a terminal vertex and which has no branch point amongst
its interior vertices. (A terminal chain need not be maximal; it can
be part of a longer terminal chain.)  A terminal edge is one meeting
a terminal vertex; it is a terminal chain of length one. (A terminal
vertex is frequently called a leaf.)

The matching number of a graph is the sum of those of its
components, so to calculate $\mn(\G)$, we may assume that $\G$ is
connected. The \emph{star} of a vertex $v$, denoted $\G_v$, is the
subgraph of $\G$ consisting of all edges meeting $v$ (and their
vertices). If $\G = \G_v$ for some $v$ then $\G$ itself will be
called a star. Its matching number is then 1. Suppose now that $\G$
is not a star. If there is a terminal chain $C$ of length 2 then a
maximum matching of $\G \setminus C$ together with a terminal edge
of $C$ gives a maximum matching of $\G$, so removing $C$ reduces the
matching number by exactly 1; similarly if there is an isolated
edge. This holds even if $\G$ has loops. However, if $\G$ is a tree
and there is neither a terminal chain of length 2 nor an isolated
edge, then there is a branch vertex $v$ such that at most one edge
meeting $v$ is not terminal. Any maximum matching of $\G \setminus
\G_v$ can be enlarged to one of $\G$ by adjoining to it any terminal
edge of $\G_v$, so removing all edges meeting $v$ also reduces the
matching number by exactly one. With this \emph{pruning procedure}
one can inductively compute the matching number of a tree, and
therefore also of a forest, whose matching number is the sum of
those of its trees. Denote the number of vertices of $\G$ by
$\operatorname{vx}(\G)$ and define the \emph{index of a graph} $\G$,
denoted $\operatorname{idx}(\G)$ to be $\vx(\G) - 2\,\mn(\G)$.
Viewing graphs as categories the definition of a product is evident
and we conjecture that $\idx(\G_1\times \G_2) = \idx(\G_1) \cdot
\idx(\G_2)$.

The \emph{adjacency matrix} of a graph $\G$ has rows and columns
indexed by the vertices of $\G$ with 1 in the $(i,j)$ place if
vertex $i$ is connected by an edge to vertex $j$. This matrix is
symmetric; its spectrum, called that of $\G$, has been intensively
studied. Suppose, however, that $\G$ is a directed graph. We then
define its \emph{skew adjacency matrix} $M(\G)$ to have $+1$ in the
$(i,j)$ place if there is an arrow directed from $i$ to $j$, to have
$-1$ in that place if there is an arrow from $j$ to $i$, with 0
there otherwise; its rank is denoted $\rk(\G)$.  This matrix defines
a skew bilinear form $B_{\G}$ on the vector space $KV$ spanned over
the field $K$ by the vertices $V$ of $\G$. Conversely, if we have a
skew bilinear form $B$ on a vector space $\mathcal{V}$ then
$\mathcal{V}$ is a direct sum of hyperbolic planes (2-dimensional
subspaces spanned by elements $v, w \in \mathcal{V}$ with $B(v,w) =
1$) and its radical (those $v$ with $B(v,w) = 0$ for all $w \in
\mathcal{V}$). Therefore skew bilinear forms on vector spaces of a
fixed dimension are completely determined up to isomorphism by their
rank, $\rk(B)$, and with suitable choice of basis are representable
by directed graphs $\G$ or their skew adjacency matrices $M(\G)$.

When $\G$ is directed we can `forget' the direction of its arrows to
get an undirected graph $|\G|$ but directed graphs $\G$ and $\G'$
with $|\G| = |\G'|$ may have different ranks. For example, if $\G$
is a square with arrows directed cyclically then $\rk(\G) = 2$ but
if  $\G'$ is obtained by reversing one arrow then $\rk(\G') = 4$.
However, for trees we have the following.

\begin{Theorem} If $\G, \G'$ are directed trees with $|\G| = |\G'|$
then their skew adjacency matrices are conjugate by a diagonal
matrix each diagonal entry of which is $\pm 1$. In particular,
$B_{\G}$ is isomorphic to $B_{\G'}$ and $\rk(\G) = \rk(\G')$.
\end{Theorem}
\noindent\textsc{Proof.\,} Suppose that $v, v'$ are vertices of a
directed tree $\G$ with an arrow $v \to v'$.  Removing this arrow
disconnects the tree. Each vertex which remains connected to $v'$
represents a basis element; changing the direction of the arrow is
essentially the same as replacing each of these basis elements by
its negative. $\Box$
\bigskip

A related concept to matching is that of a \emph{node cover}, i.e.,
a set $T$ of vertices of $\G$ such that every edge has an end in
$T$. Here we wish to minimize $T$. Minimal node covers may have
different sizes; one which achieves the absolute minimum is a
\emph{minimum cover} and its size, the \emph{cover number}, will be
denoted $\cn(\G)$. If $M$ is a matching and $T$ a node cover, then
no vertex in $T$ can cover more than one edge in $M$. Therefore $\#M
\le \#T$, so $\mn(\G) \le \cn(\G)$. In general this inequality is
strict. (The triangle has vertex number equal to $2$ but matching
number $1$.) However, for a \emph{bipartite} graph, i.e., one in
which the vertices can be partitioned into two disjoint sets with
every edge connecting a vertex in one to a vertex in the other,
Menger's Theorem asserts that equality holds. One may think of a
bipartite graph as having its vertices colored, say either black or
white, with edges always connecting vertices of different colors.
Trees are bipartite.

\section{The graph $\G(S)$}\label{sec:graph}
Suppose that $\g$ is a Lie subalgebra of $\spl(n)$ which contains
the Cartan subalgebra $\h$ of traceless diagonal matrices. It is
then spanned by these matrices together with all $e_{ij}\in\g$. Let
$\Pi(\g)$ denote the support of $\g$, i.e., the set of those pairs
of indices $(i,j)$ for which $e_{ij} \in \g$, together with the
diagonal pairs $(i,i)$; for $\g = \PC(n,m)$ we will write simply
$\Pi(n,m)$. Given a small functional $F_S$ on $\g$ (where tacitly $S
\subset \Pi(\g)$) we now define a directed graph $\G(S)$ with
vertices the union of (i) the set of those $e_{ij}$ in $\g$ where
$(i,j)$ is neither in $S$ nor on the diagonal and (ii) a vertex
labeled $d_s$ for each $s\in S$. Writing $B_S$ for $B_{F_S}$ draw an
arrow $e_{ij} \to e_{kl}$ whenever $B_S(e_{ij}, e_{kl}) = 1$, i.e.,
if $j=k$ and $(i,l) \in S$, and for every $s\in S$ add an arrow
$d_s\to e_s$ since we also have $B_S(d_s, e_s) = 1$.

The decomposition of $\g$ (which here need not be Frobenius) into
eigenspaces of $\ad(D_S)$ also decomposes $\G(S)$ into disjoint
subgraphs, one for each pair of eigenvalues $(m, 1-m)$, since an
$e_{ij}$ in the eigenspace for the eigenvalue $m$ can only be linked
to one for the eigenvalue $1-m$ while $d_s$ can only be linked to
$e_s$ (the eigenvalues for which are $0$ and $1$, respectively). The
eigenspace components of $\G(S)$ may themselves decompose further.
Each component and hence all of $\G(S)$ is therefore bipartite. When
$\g$ is Frobenius one can in principle use $\G(S)$ to calculate the
associated solutions to the CYBE, and we do so in several cases. The
graph $\G(S)$ then effectively organizes the inversion of $B_S$,
which may be large but sparse. With the preceding notations we have
the following.

\begin{Theorem} \label{th:rk} If $\G = \G(S)$ is a tree then
$\rk(\G) = \rk(B_S) = 2\mn(|\G|)$. If moreover $F_S$ is regular then
$\idx(\G_S) = \idx(\g).$
\end{Theorem}
\noindent\textsc{Proof.\,} Note that while $\G$ is a directed graph,
by the remarks in the preceding section its rank does not depend on
the directions of the arrows so it must depend only on $|\G|$, so
when there can be no confusion we will write simply $\G$ for $|\G|$.
The theorem is obvious for $\G$ a star. (Taking the vertex at the
center of the star as first basis element the matrix of $B$ has
non-zero elements only in the first row and first column.)  For
larger $\G$ suppose that the theorem holds for all trees with fewer
vertices than $\G$ and apply the pruning procedure of the previous
section. If there is a terminal chain $C$ of length 2 then by
hypothesis the theorem holds for $\G\setminus C$, the matching
number of which is one less than that of $\G$. Take the terminal
vertex which has been removed and that connected to it (also
removed) as the first and second basis elements. Since only the
second is linked to any elements of $\G\setminus C$ it is evident
that the matrix corresponding to $\G\setminus C$, which is obtained
from that of $\G$ by removing the first and second rows and columns,
has rank exactly two less than that of $\G$. When the pruning
removes a terminal star $\G_v$ the argument is a slight elaboration
of the preceding. Take $v$ as the first vertex, followed by the
others of the star and then those of $\G\setminus\G_v$. Since the
only vertex of $\G_v$ linked to $\G\setminus\G_v$ is $v$, removing
the first row and column of the full matrix of the form reduces it
to the direct sum of a zero matrix and that corresponding to
$\G\setminus \G_v$. It is evident then that the difference in ranks
is again exactly two, regardless of the number of vertices in
$\G_v$. This proves the first assertion of the theorem; the second
follows. $\Box$
\bigskip

The graph $\gamma(S)$ is always a tree and we should like to be able
to show that this implies the same for the components of $\G(S)$
(which are, in any case, bipartite), but that for the moment is an
open question. The case of most interest for us is, however, that
where $\gamma(S)$ is a rooted tree (all arrows being directed away
from the root).

\begin{Theorem} If $\gamma(S)$ is a rooted tree then so are all
components of $\G(S)$.
\end{Theorem}
\noindent\textsc{Proof.\,} Observe first that if $s\in S$ then the
only arrow that can terminate on $e_s$ is $d_s\to e_s$, so in the
proof we may disregard all $d_s$ (which will be the root of any
component in which it appears) and consider only vertices $e_{ij}$
of $\G(S)$ with $i\ne j$. Suppose, if possible, that we had arrows
$\e_{ij} \to e_{jk} \ot e_{hj}$ in $\G(S)$. Then in $\gamma(S)$ we
must have $i\to k \ot h$, contradicting the assumption that
$\gamma(S)$ is a rooted tree with all arrows directed away from the
root.$\Box$
\bigskip

\section{Isotropic and Lagrangian subspaces}\label{sec:Lagrange}
Let $V$ denote a vector space over a field of characteristic
different from $2$ endowed with a skew bilinear form $B(-,-)$. A
subspace $W$ is called isotropic if $B(w_1,w_2)=0$ for all $w_1,w_2
\in W$; equivalently it is contained in its orthogonal complement
$W^{\perp}= \{w'\in V|B(w,w') = 0, \text{all } w\in W\}$. If $B$ is
non-degenerate then $\dim V = 2\ell$ is even and every maximal
isotropic subspace $L$ has dimension $\ell$; such subspaces are
called Lagrangian. One has $V = L\oplus L^{\perp}$, where
$L^{\perp}$ is again Lagrangian and choosing any basis
$\{v_1,\dots,v_{\ell}\}$ in $V$ there is a dual basis
$\{v_1',\dots,v_{\ell}'\}$ with $B(v_i,v_j') = \delta_{ij}$. The
matrix of $B$ relative to the basis
$\{v_1,\dots,v_{\ell},v_1',\dots,v_{\ell}'\}$ then has the form
$\begin{pmatrix} 0 & I_{\ell}\\-I_{\ell}&0\end{pmatrix}$ the inverse
of which is just its negative.

Suppose that we have a Frobenius Lie algebra $\f$ with Frobenius
functional $F$ and form $B = B_F$. The computation of the associated
$r$-matrix becomes trivial if we can write $\f$ as a direct sum of a
pair of Lagrangian subspaces $\f = L \oplus L'$ with an explicit
duality between them. For with a basis of $\f$ consisting of  a
basis $x_i$ of $L$ followed by the dual basis $x_i'$ of $L'$ the
matrix of $B$ has the form $\begin{pmatrix} 0 &
I\\-I&0\end{pmatrix}$, where $I$ is an identity matrix of size the
dimension of $L$. The inverse of this matrix is just its negative,
and since we could multiply by any non-zero scalar, we can simply
take $r = \sum x_i\wedge x_i'$.

\begin{Theorem} A Frobenius Lie algebra $\f$ whose principal element has only
integral eigenvalues can be decomposed into a direct sum of
Lagrangian subspaces in duality with each other, one of which is a
Lie subalgebra and the second a module over the first under the Lie
multiplication. In particular this holds for Frobenius subalgebras
of $\spl(n)$ which contain the Cartan subalgebra.
\end{Theorem}

\noindent\textsc{Proof.\,} It follows from \cite{G:principal} (or
explicit trivial calculation) that principal elements $\hat F$ have
adjoints $\ad\hat F$ whose eigenvalues are integers. Letting $\f_m$
denote the eigenspace for the integer $m$, the dual of which is
$\f_{1-m}$, one see that $\fe = \sum_{m \mathrm{\, even}}\f_m$ and
$\fo = \sum_{m \mathrm{\, odd}}\f_m$ are dual Lagrangian subspaces.
Since $[\f_m, \f_n] \subset \f_{m+n}$, $\fe$ is a Lie subalgebra and
$\fo$ is a module over $\fe$ under the multiplication in $\f$.
$\Box$
\bigskip

Suppose that one can readily find only one Lagrangian subspace $L$
(not necessarily a subalgebra) with a complement $L'$ which is not
necessarily Lagrangian. If a basis $\{v_1,\dots,v_{\ell}\}$ of $L$
is given then we can still find in $L'$ a dual basis
$\{v_1',\dots,v_{\ell}'\}$ but now the matrix of $B$ relative to the
basis $\{v_1,\dots,v_{\ell},v_1',\dots,v_{\ell}'\}$ will have the
form $M=\begin{pmatrix} 0 & I_{\ell}\\-I_{\ell}&Q\end{pmatrix}$ for
some $\ell \times \ell$ matrix $Q$. Then $M^{-1}=
\begin{pmatrix} Q & -I_{\ell}\\I_{\ell}&0\end{pmatrix}$ so no
computation is required to find the associated $r$-matrix beyond
finding the dual basis to $\{v_1,\dots,v_{\ell}\}$.

\section{The cyclic functional}\label{sec:cyclic} The
\textit{cyclic functional} on $\PC(n,m)$ (with $(n,m) = 1$) is the
one in which we are most interested.  Its associated tree graph will
be denoted simply by $\gamma = \gamma(n,m)$. To define it, write the
integers $1,\dots n$ in the order $1, m+1, 2m+1, \dots,
\linebreak[3](n-1)m+1, 1$ where all $km+1$ are understood to be
reduced modulo $m$ except $n$ itself, and the order is understood to
be cyclic. For example, with $n=13, m=5$ one has $\dots
1,6,11,3,8,13,5,10,2,7,12,4,9,1,6,11,\dots$. As in
\cite{GG:Boundary}, divide the array into \textit{strings} of
consecutive entries inside which the integers are increasing in
their natural order. In the example the strings are $(1,6,11),
(3,8,13), (5,10), (2,7,12), (4,9)$; this will be called the first
\textit{cycle of strings}. The directed graph $\gamma$ as before has
vertices the integers $1,\dots,n$ to which we add directed edges
beginning with those joining successive integers in each string. In
the illustration we add $1\to6\to11,\, 3\to8\to13,\,5\to10,\,
2\to7\to12$. Now consider the ``drops" in the original cyclic order,
i.e., those integers which are less than their predecessors (the
first integers of the strings); here they are $3,5,2,4,1$. Arrange
these similarly in ascending strings, $(1,3,5),(2,4)$ (note the
cyclic order), giving the second cycle of strings. Now draw arrows
from the \emph{larger} to the \emph{smaller} integers in the
strings: $1\ot 3\ot 5,\, 2 \ot 4$. Next, instead of removing the
drops, remove the ``rises", those numbers which, reading in
\textit{reverse} cyclic order, are larger than their predecessors
(the last elements in the strings). Here one has only $5, 4$. Now go
back to the first procedure, grouping them into ascending strings.
In this case, we have only one arrow to add, $4\to 5$. (Remember
that the cyclic order is $4,5,4$.) In general we continue
alternating between drawing arrows from smaller integers to larger
ones and from larger integers to smaller ones until the process ends
with but a single group. This also defines the set $S = S_{n,m}$
supporting the cyclic functional, denoted $F_{n,m}$; it is the set
of those pairs $(i,j)$ for which there is an arrow from $i$ to $j$
in $\gamma$. Isomorphic algebras, in particular $\PC(n,m)$ and
$\PC(n,n-m)$, generally do not have isomorphic trees. The simplest
example of this is $\PC(5,2)$, whose tree has root $2$, at which it
branches, and $\PC(5,3)$, whose tree has root $3$ but there is no
branching there. Note that the $m$th superdiagonal is filled by $S$,
i.e., all $(i, i+m)$ in $\Pi(n,m)$ are in $S$;

\begin{Theorem}The graph $\gamma(n,m)$ is a
rooted tree with all branches directed away from the root; with the
exception of the unique root, every $j \in \{1,\dots,n\}$ has a
unique immediate predecessor $i$ such that $i\to j$. \end{Theorem}

\noindent\textsc{Proof.\,}An integer will appear in a cycle of
strings only if it has no predecessor in any preceding cycle, so for
any integer there can be at most one predecessor. Since the process
continues until the cycle consists of just a single string, there
can only be one integer in $\{1,\dots, n\}$ without a predecessor,
and that will be the root. $\Box$
\bigskip

We will say more generally that a vertex $v$ of a directed graph has
another vertex $w$ as  a predecessor if there is an arrow $w\to v$.
Viewing $S$ for the moment as an $n\times n$ matrix, the theorem
shows that there is precisely one non-zero element (equal to 1) in
each column except for that column corresponding to the root (which
in the example is column 4). It follows for $\G(n,m)$ that
predecessors, when they exist, are unique, so each component of
$\G(n,m)$ is therefore itself a rooted tree. This is an essential
feature of the cyclic functional.

The directed tree $\gamma = \gamma(n,m)$ of the cyclic functional
can be built recursively. If $n > 2m$ then $\gamma(n-m,m)$ is a
subgraph of $\gamma(n,m)$ and the latter is obtained from the former
by attaching an outgoing arrow to every vertex $i$ with
\linebreak[4] $n-2m+1 \le i \le n-m$ with the integer $i+m$ at the
point of the arrow. Note that $n-m+1,\dots,m$ are the ends of
$\gamma(n,m)$ (vertices with no outgoing arrow) so one may view the
reduction as the removal of all ends together with the unique arrows
terminating on them. In this case, which we will call \emph{stable
reduction} the root of $\gamma(n-m,m)$ coincides with that of
$\gamma(n,m)$. If $n-m < m$ then stripping all the $n-m$ ends (whose
labels are $m+1,\dots,n$) from $\gamma(n,m)$ produces a graph which
as a directed tree but without labeled vertices is identical with
$\gamma(m, 2m-n)$. To get the correct labeling one must now replace
every $i$ with $m+1-i$. It follows that if $\rho(n,m)$ is the root
of $\gamma(n,m)$ and $\rho(m, 2m-n)$ that of $\gamma(m, 2m-n)$ then
$\rho(n,m) = m+1-\rho(m, 2m-n)$. Since the numbering and root both
change, we call this case \emph{unstable reduction}. For example,
$\gamma(17,6)$ reduces stably to $\gamma(11,6)$, so the root does
not change, but $\gamma(11,6)$ reduces unstably to $\gamma(6,1)$.
The root of the latter is $1$ (and that of any $\gamma(n,1)$ is
clearly $1$), so the root of $\gamma(11,6)$ is $(6+1)-1 = 6$, which
therefore is also the root of $\gamma(17,6)$.

\begin{Theorem} The cyclic functional is Frobenius. \end{Theorem}
\textsc{Proof.\,} We show that $F_{n,m}$ is Frobenius if and only if
the same is true for $F_{n',m'}$, where $(n',m') = (n-m,m)$ if $n >
2m$ (stable case) or $(n',m') = (m,2m-n)$ if $n<2m$ (unstable case);
ultimately $(n,m) = (2,1)$ where the theorem is obvious. The
reduction processes are similar in that in each case we remove two
disjoint blocks with the same numbers of elements from $\Pi(n,m)$.
In the stable case, the first block to be removed consists of the
last $m$ rows of $\Pi(n,m)$; it is a block of $m$ rows and $n-m$
columns. After that the second block to be removed consists of the
last $m$ columns of what remains; it is a block of $n-m$ rows and
$m$ columns. Note that the first block contains no element of $S$,
while the second contains one element of $S$ in each column. The
resulting configuration is that of $\Pi(n-m,m)$ in its standard
position. In removing the second block we removed $m$ elements of
$S$. Those that remain are in the correct positions for $S_{n-m,m}$.
In the unstable case, the first block to be removed consists of
those $(i,j)$ with $n-m+1 \le i \le n$, $m=1 \le j \le n$; it is a
block of $m$ rows and $n-m$ columns. After that the second block
removed consists of those $(i,j)$ with $1 \le i \le n-m, n-m+1 \le j
\le n$; it is a block of $n-m$ rows and $m$ columns. Again, the
first block contains no element of $S$ but now the second block
contains one element of $S$ in each row. The configuration that
remains in the unstable case is not that of any usual $\Pi(n', m')$,
but after rotation through a half circle becomes that of $\Pi(m,
n-m)$. Since there are now $m$ rows and columns, this amounts to
replacing $(i,j)$ in the reduced configuration by $(m+1-i, m+1-j)$.
In removing the second block we removed $n-m$ elements of $S$. After
the rotation, they are in the correct position for $S(m, 2m-n)$.
Note that in both cases the two blocks have the same number of
elements, the first block contains no elements of $S$, the second
block contains no elements on the diagonal, and that the number of
diagonal elements of the first block is equal to the number of
elements of $S$ in the second.

These steps are illustrated in the following figure, where the first
block is marked by the solid black squares $\blacksquare$ and the
second by the open ones $\square$. Entries in $S$ are marked by
$\mathbf{x}$, including those inside the open squares. (There can be
none in the locations marked by the black squares.) Entries outside
$\Pi(7,3)$  must be zero and are marked by dots.

\begin{equation}\label{eqn:redux}
{\begin{pmatrix}
\ast&\ast&\ast&\mathbf{x}& \square&\square&\square\\
\ast&\ast&\mathbf{x}&\ast&\boxtimes&\square&\square\\
\mathbf{x}&\ast&\ast&\ast&\square&\boxtimes&\square\\
\cdot&\cdot&\cdot&\ast&\square&\square&\boxtimes \\
\cdot&\cdot&\cdot&\blacksquare&\blacksquare&\blacksquare&\blacksquare\\
\cdot&\cdot&\cdot&\blacksquare&\blacksquare&\blacksquare&\blacksquare \\
\cdot&\cdot&\cdot&\blacksquare&\blacksquare&\blacksquare&\blacksquare
\end{pmatrix}}\rightsquigarrow {\begin{pmatrix}
\ast&\square&\square&\boxtimes\\
\ast&\ast&\mathbf{x}&\blacksquare\\
\mathbf{x}&\ast&\ast&\blacksquare\\
\cdot&\cdot&\cdot&\blacksquare
\end{pmatrix}}\rightsquigarrow
\begin{pmatrix}\ast&\cdot&\cdot\\
\ast & \ast &\mathbf{x} \\
\mathbf{x} & \ast & \ast
\end{pmatrix}
\end{equation}
\textsc{Figure 1.}\texttt{\ \small{Steps in the reduction of
$\Pi(7,3)$: Stable reduction to $\Pi(4,3)$, followed by unstable
reduction to $\Pi(3,2)$. The last matrix must be rotated through
$180^{\circ}$ to bring $S$ into the standard position for
$\gamma(3,2)$.}}
\bigskip

Consider next what happens to the graph $\G(S)$ in the reduction
process. The first (black square) block removed contains no element
of $S$, so none of its entries has a successor (an element of which
it is the predecessor), but every element of the first block which
is not on the diagonal has a (necessarily unique) predecessor in the
second (white square) block, and that predecessor can not be an
element of $S$. This defines a bijective map from the non-diagonal
elements of the first block to those elements of the second block
which are not in $S$. An element of the second block may have no
predecessor, but if it does, that predecessor (which may be an
element of $S$) is not contained in the second block and therefore
remains after the reduction. Letting $\G' = \G(S')$ denote the graph
of what remains after the reduction, to build $\G = \G(S)$ from
$\G'$ one must do the following: First, if $(k,l)$ is an element of
the first (black square) block whose predecessor $(j,k)$ in the
second (white square) block itself has a predecessor $(i,j)$ then
one must attach a terminal chain of length $2$: $e_{ij}\to e_{jk}\to
e_{kl}$. (Here $e_{ij}$ is necessarily a vertex of $\G(S')$.) If
$(k,l)$ is an element of the first  block whose predecessor $(j,k)$
in the second block has no predecessor then adjoin a disconnected
arrow $e_{jk}\to e_{kl}$. Finally, for each $s \in S\setminus S'$
adjoin another disconnected arrow from a new vertex labeled $d_s$
(the diagonal element dual to $e_s$) to $e_s$. Note that in every
case when an arrow points from a vertex $v$ to its successor $v'$
one has $F_S([e_v, e_{v'}]) = 1$, and conversely.

It is evident now that the matching number of $\G$ is the matching
number of $\G'$ plus half the number of elements in the two blocks
removed (or the number in either single one of them) since to every
pair of elements (non-diagonal element of first block, its
predecessor) or (diagonal element of first block, element of $S$ in
its column) we have either attached to some element of $\G'$ a
terminal chain of length $2$ or have adjoined a disjoint link. The
reduction brings us to some smaller $\Pi(n',m')$ (where again
$(n',m') = 1$) with its corresponding cyclic form; if the latter is
Frobenius then so is that with which we started. Since the smallest
case is evident on inspection, this ends the proof. $\Box$
\bigskip

In the proof we have, in effect, recursively constructed the $\G$
there from $\G'$. We conclude this section by recursively
constructing $\G(7,3)$ starting with $\G(3,2)$: $\{e_{12} \to
e_{23}, d_{13}\to e_{13},  d_{21}\to e_{21}\}$; here $d_{13} = \e_1
+\e_2$ and $d_{21} = \e_2$. Since the reduction from $\G(4,3)$ to
$\G(3,2)$ was unstable we must now complement all indices with
respect to $4$ ($d_{13}$ becomes $d_{31} =\e_3+ \e_2$ and $d_{21}$
becomes $d_{23} = \e_2$) and then adjoin the new arrows. To indicate
which arrows are new and which have come from (the rotated)
$\G(3,2)$ we will indicate the latter by a double arrow
$\Rightarrow$ and the new arrows by a single arrow $\rightarrow$.

$$e_{32}\Rightarrow e_{21}\rightarrow
e_{13}\rightarrow e_{34}\quad d_{23}\Rightarrow e_{23}\quad
d_{31}\Rightarrow e_{31}\quad d_{14} \to e_{14}\quad
e_{12}\rightarrow e_{24}$$
\begin{center}\textsc{Figure 2.}\texttt{\ The graph
$\G(4,3)$}\end{center}
\medskip

 Note that
$d_{31}$ is still $e_3+ \e_2$ and $d_{23}$ remains $\e_2$, while
$d_{14} = \e_1+\e_2+\e_3$. There are $12$ vertices (the dimension of
$\PC(4,3)$) and the matching number of this graph is $6$. The graph
below of $\G(7,3)$ is arranged so that arrows stemming from
$\G(3,4)$ are vertical and new ones, denoted by $\dashrightarrow$,
are horizontal. There are $11$ components, $36$ vertices (dimension
of $\PC(7,3)$), and the matching number can be seen to be $18$.
\begin{gather*}
\begin{array}{ccccc}
e_{32}&\dashrightarrow& e_{26}&\dashrightarrow &e_{65}\\
\Downarrow&&&&\\
e_{21} &\dashrightarrow &e_{15} &\dashrightarrow &e_{54}\\
\downarrow&&&&\\ e_{13}&&&&\\ \downarrow&&&&\\ e_{34}&
\dashrightarrow &e_{46}& \dashrightarrow &e_{67}
\end{array}\qquad
\begin{array}{ccccccc}
d_{23}&&&&\\ \Downarrow &&&&\\
e_{23}&\dashrightarrow&\e_{35} &\dashrightarrow&e_{56}\\ &&&&&\\d_{31}&&&&\\ \Downarrow&&&&\\
e_{31}&\dashrightarrow& e_{16} &\dashrightarrow&e_{64}\\&&&&&\\e_{12}&&&&\\
\downarrow &&&&
\\e_{24}&\dashrightarrow&\e_{45}
&\dashrightarrow&e_{57}\\&&&&&
\end{array}\qquad
\begin{array}{c}
d_{14}\\ \downarrow \\e_{14}
\end{array}\\
d_{25}\dashrightarrow e_{25}\qquad d_{36}\dashrightarrow
e_{36}\qquad
d_{47}\dashrightarrow e_{47}\\
e_{17}\dashrightarrow e_{74}\qquad e_{27}\dashrightarrow e_{75}
\qquad e_{37} \dashrightarrow e_{76}
\end{gather*}
\begin{center}\textsc{Figure 3.}\texttt{\ The graph $\G(7,3)$}
\end{center}
\bigskip
The eigenvalue pair $(m,1-m)$ to which each component of this graph
belongs is determined by its root; the values are easily obtained by
using the principal element $D_{S_{\mathrm{cyclic}}}$ of equation
(\ref{eqn:D_S}). (For example, the root $e_{32}$ of the largest
component belongs to the eigenspace with $m=-1$.) The $\G(S)$ here
has a unique perfect matching; this will always be the case for any
small Frobenius functional $F_S$.

\section{Solutions to the CYBE from graphs}\label{sec:CYBEgraph}
The solution to the CYBE derived from $\PC(n,m)$ using the cyclic
functional will be denoted $r(n,m)$. It is naturally a sum of terms
corresponding to the components of $\G(n,m)$. The bilinear form
defined by the functional will for the moment be denoted simply $B$.
\subsection{$\PC(4,3)$ and $\PC(7,3)$} The graph $\G(4,3)$  (Figure
2) has five components of which the last four are just isolated
links. The four corresponding summands of the solution to the CYBE
are  $d_{23}\wedge e_{23},\ d_{31}\wedge e_{31},\ d_{14} \wedge
e_{14}, \ e_{12}\wedge e_{24}$. To compute the summand for the first
component, which illustrates the more general situation, notice that
a terminal chain of length 2 has been added to a component of the
graph $\G(3,2)$; the reduction process of the preceding section
shows that any $\G(n,m)$ (where $(n,m) = 1$) is built from a smaller
one by the addition of isolated links and terminal chains of
length 2. The graph encodes that $B(e_{13},e_{34})=1$ and
$B(e_{21},e_{13}) = 1$, so  $B(e_{21}+e_{34}, e_{13})=0$. Since
$B(e_{21},e_{34}) = 0$ we also have $B(e_{32},e_{21}+e_{34}) = 1$.
We can therefore separate the first component into two components,
$e_{32}\Rightarrow e_{21}+e_{34}$ and $e_{13}\to e_{34}$. By
replacing $e_{21}$ (to which the terminal chain was attached)  by
$e_{21}+e_{34}$ we have in effect detached the terminal chain (which
is now reduced to an isolated link) from the previous graph
$\G(3,2)$. What remains is now also just an isolated link so the
more general process of detaching terminal chains of length two by
which the larger graph has been built from the smaller is ended. The
summand corresponding to the first component of $\G(3,4)$ is thus
$e_{32}\wedge(e_{21}+e_{34}) + e_{13}\wedge e_{34}$. Inserting the
values for the various $d_s$ we therefore have
\begin{equation*}
\begin{split}
r(4,3)\  = &\,e_{32}\wedge(e_{21}+e_{34}) + e_{13}\wedge e_{34} +\\
&\ \ \ \e_2\wedge e_{23}+ (\e_2+\e_3)\wedge e_{31}+ (\e_1+\e_2+\e_3)
\wedge e_{14} +e_{12}\wedge e_{24}.
\end{split}
\end{equation*}
Analyzing the first component of $\G(7,3)$ in the same way,
disconnecting the three terminal chains (dotted arrows) that have
been added to a component of $\G(4,3)$ reduces it to
\begin{gather*}
(e_{32}+e_{65})\Rightarrow(e_{21}+e_{54})\to e_{13} \to
(e_{34}+e_{67})\\e_{26}\dashrightarrow e_{65}\qquad
e_{15}\dashrightarrow  e_{54} \qquad e_{46}\dashrightarrow e_{67}
\end{gather*}
We can now treat the first line above exactly as before, so the
contribution of the first component of $\G(7,3)$ to $r(7,3)$ is
$$(e_{32}+e_{65})\wedge [(e_{21}+e_{54})+(e_{34}+e_{67})]
+e_{13}\wedge (e_{34}+e_{67}) +e_{26}\wedge e_{65}+ e_{15}\wedge
e_{54} + e_{46}\wedge e_{67}$$ It is not difficult now to compute
the 13 additional summands of $r(7,3)$; we omit it. (The isolated
links contribute 7 and the remaining three components each
contribute two.)

\subsection{$\PC(n,1)$}  The cyclic functional $F$ for $\PC(n,1)$
reduces to what is sometimes called the \emph{prime functional}. Its
carrier is $S = \{(1,2), (2,3),\dots,(n-1,n)\}$. One has
\begin{gather*}
\gamma(S) = 1\to2\to\dots\to (n-1)\to n, \quad  d_{(i,i+1)} = \e_1 +
\e_2 +\dots + \e_i, \\  \begin{align*}D_S &= (n-1)\e_1 + (n-2)\e_2 +
\dots +\e_{n-1} \\&=(1/2)[(n-1)e_{11} +(n-3)e_{22} + \dots
+(3-n)e_{n-1,n-1} +(1-n)e_{nn}].\end{align*}
\end{gather*}
 The eigenspace of $\ad(D_S)$ for
the eigenvalue $m$ is the Cartan subalgebra when $m=0$ and otherwise
the $m$th superdiagonal (subdiagonal if $m$ is negative). A closed
form for the associated solution to the CYBE was given in
\cite{GG:Boundary}. Using the graph $\G =\G(n,1)$ we see here why it
has the given form. In this simple case $\G$ is a disjoint union of
chains.
\begin{gather*}
d_{i,i+1} \to e_{i,i+1}, \qquad i = 1,\dots,n-1\\
(1,3)\to(3,2)\to(2,4)\to(4,3)\to\dots\to(n-2,n)\to (n,n-1)\\
(1,4)\to(4,2)\to(2,5)\to(5,3)\to\dots\to (n-3,n)\to (n,n-2)\\ \dots\\
(1,n-1)\to(n-1,2)\to(2,n)\to (n,3)\\(1,n)\to (n,2)
\end{gather*}
These are, respectively, the components of $\G$ for the eigenvalue
pairs \linebreak $(0,1)$ (the $n-1$ short chains of the first row)
$(2,-1), (3,-2),\dots,(n-2,3-n), (n-1, 2-n)$. Applying the procedure
described in the preceding section, from the last chain one has
$e_{1n}' = e_{n2}$, from the next that $e_{2n}' = e_{n3}$ and
$e_{1,n-1}'=e_{n-1,2} + e_{n,3}$, and so forth. From the first chain
we get, in particular, that $e_{1,3}' = e_{3,2} + e_{4,3} + \cdots
+e_{n,n-1}$. Collecting terms $\sum_{i<j}e_{ij}\wedge e_{ij}'$ and
noting that $d_{i,i+1}\to e_{i,i+1}$ (so $d_{i,i+1}'= e_{i,i+1}$) we
get the closed form of \cite{GG:Boundary} for the solution to the
CYBE associated to $\PC(n,1)$. As given there, with terms collected
in slightly different order and writing simply $d_p$ for $d_{p,p+1}$
it is
$$r(n,1)\quad = \quad
\sum_{p=1}^{n-1}d_p\wedge e_{p,p+1} +\sum_{i<j}
\sum_{m=1}^{j-i-1}e_{i,j-m+1}\wedge e_{j,i+m}.$$ The carrier $S$ of
the prime functional is unchanged by reflection across the
antidiagonal, so the prime functional will work also for
$\PC(n,n-1)$. However, $m=1$ and $m=n-1$ are the only cases where
the prime functional is Frobenius; in other cases, while the
eigenspace for $n-1$ always has dimension $1$ that for  $2-n$ (which
should be its dual) vanishes.

\subsection{$\PC(n,2)$} Note that $n$ must be odd. Here $\G(n,2)$
again consists only of
chains. One has
\begin{gather*}S=\{(2,1),(1,3),(2,4),(3,5),\dots,(n-2,n)\}, \\
\gamma(2,n)\quad =\quad (n-2)\ot (n-4)\ot \cdots \ot 4 \ot 2
\to1\to3 \to5\cdots\to n.
\end{gather*}
From this one finds that
\begin{align*} d_{21} &=
\e_2+\e_4+\dots+\e_{n-2}, \\d_{2k+2,\,2k+3}& =
-(\e_{2k+2}+\e_{2k+4}+\cdots+\e_{n-3} +\e_{n-1}),\\
  d_{2k+1,\,2k+3}&= \e_1+\e_2+\e_3+\cdots +e_{2k+1}-d_{2k+2,2k+3}\\&=
(\e_1+\e_2+\e_3+\cdots
+e_{2k+1})+(\e_{2k+2}+\e_{2k+4}+\cdots+\e_{n-3}
+\e_{n-1})\end{align*}

For $\G(n,2),\, n\ge5$ one has
\begin{gather*}
d_s\to e_s,\quad \mathrm{all} \quad s \ne (2,1)\\
e_{12}\to e_{23}\to  e_{34} \to \cdots \to e_{n-1,n}\\
d_{21}\to e_{21}\to e_{14}\to e_{43} \to \cdots \to e_{n-4,n-1} \to
e_{n-1,n-2}\\
\mathrm{and\ for\ all\ }5\le j\le n ,\,k\ge0,\\
e_{1,j}\to e_{j,3} \to e_{3,j+2}\to e_{j+2,5}\to \cdots \to
e_{2k+1,2k+j} \to e_{2k+j, 2k+3}\to \cdots  \\
e_{2,j}\to e_{j,4} \to e_{4,j+2}\to e_{j+2,6}\to \cdots  \to
e_{2k+2,2k+j} \to e_{2k+j, 2k+4}\to \cdots
\end{gather*}
where a chain terminates when any index exceeds $n$. It follows that
\begin{align*}
r(1,2) = \sum_{s\in S,s\ne(2,1)}d_s\wedge e_s + &d_{21}\wedge
(e_{21} +\sum_{k\ge2}e_{2k,2k-1})\\& +\sum_{i<j,j\ne
i+2}e_{ij}\wedge \sum_{k\ge0}e_{j+2k,i+2k+2}
\end{align*}
where the sums terminate when the indices are out of range. The
principal element is
$$\hat F =
\operatorname{diag}(0,1,-1,0,-2,-1,-3,\dots,(5-n)/2, (1-n)/2) +
[(n-1)(n-3)/4n]\,\mathbf{I_n}.$$
\section{Some other small Frobenius functionals}\label{sec:functionals} There are many small
Frobenius functionals that can be defined for $\PC(n,m)$. We give
several in this section just to show the variety but consider the
Dergachev-Kirillov functional separately.
\subsection{The subprime functional} For $m>1$ the
simplest functional after the prime is the \emph{subprime}, defined
by taking $S$ to be the union of the sets $(i,i+m), i = 1,\dots,n-m$
and $(i+1,i), i = 1,\dots,m-1$. No column contains more than one
element of $S$ and the only one without one is the $m$th. For $n
\equiv -1\mod m$ this is just the cyclic functional. The subprime
functional is still Frobenius when $n \equiv 1\mod m$; the proof is
by a reduction process almost identical with that in the cyclic case
and the associated $r$ matrix is still relatively simple to
construct, at least for small $n$. However it is only for $n \equiv
\pm 1\mod m$ that the subprime functional is Frobenius.

\subsection{The upper triangular functional} The
\emph{upper triangular} $S$ is a modification of the cyclic
functional with the property that $S$ is contained entirely above
the diagonal. As with the cyclic functional, first write the
integers $1,2,\dots,n$ in the order $1, m+1, 2m+1, \dots$ where the
entries are understood modulo $m$ except for $m$ itself and the
order is cyclic (so we end again with 1). For example, with
$n=12,m=5$ we have $1,6,11,4,9,2,7,12,5,10,3,8,1$ where we have
repeated the first integer to emphasize the cyclic order. If $n
> 2m$ read the sequence forwards, remove the $m$ largest numbers
from the sequence, and if $j$ is one of these and $i$ its
predecessor in the sequence, take the pair $(i,j)$ into $S$. View
what remains of the sequence as associated to $\PC(n',m')$ with $n'
= n-m, m' = m$. However, if $n < 2m$ remove the $n-m$ smallest
numbers and if $j'$ is one of these and $i'$ its predecessor take
$(j',i')$ (note the transposition) into $S$. Think of reading the
sequence backwards. Now view what remains of the sequence as
associated to $\PC(n',m')$ with $n'= m, m' = n-m$ \emph{except that
the indexing has been shifted up by} $n-m$ since $1,\dots,n-m$ have
been removed. In either case, if now $n'
> 2m'$ continue reading the sequence in the same direction in which
it was last read, otherwise reverse direction. So for $\PC(12,5)$ we
first strike $8,9,10,11,12$, taking $(3,8), (4,9), (5,10), (6, 11),
(7,12)$ into $S$ and leaving the sequence $1,6,4,2,7,5,3,1$. Now we
proceed as if we were in the case $n=7, m=5$, so we must strike the
the two smallest integers, 1 and 2, whose predecessors are
respectively 3 (remember the cyclic order!) and 2, so we take
$(1,3), (2,4)$ into $S$. The remaining sequence is $6,4,7,5,3,6$ and
we proceed as if it were the case of $n=5, m=2$. \emph{We now
continue reading the sequence in the same direction as before},
namely backwards, removing the smallest entries 3 and 4, and taking
$(3,5), (4,6)$ into $S$.  This leaves $6,7,5,6$ with $n=3, m=2$. We
must reverse direction again, and since we are again reading in the
\emph{forward} direction we strike only the one \emph{largest}
integer, 7, and take $(6,7)$ into $S$. The sequence has been reduced
to $6,5,(6)$ and we are in final case of $n=2, m=1$. We do not
reverse direction and take $(5,6)$ into $S$, producing an upper
triangular $S$ for which $\gamma(S)$ is obviously a tree. Note that
we started with $m=5$ and lastly took in $(5,6)$; the final entry
into $S$ will always be $(m,m+1)$.

The foregoing is actually a prescription for a reduction procedure
removing pairs of blocks analogous to that in the cyclic case; here,
too, it gives an inductive proof that $S$ is Frobenius. This is
illustrated in the following diagram for $\PC(12,5)$ where dots
indicate the places that must be filled with zeros and positions
which have some manner of ``x'' in them are those in $S$.  The first
entries taken into $S$ were (3,8), (4,9), (5,10), (6,11), and
(7,12); these are in the positions indicated by $\boxtimes$. The
last five rows are removed, which removes the $5\times10$ block
indicated by the black squares $\blacksquare$ and then the last 5
columns, which removes the $10\times 5$ block indicated by the
squares $\sq$ (including those in the positions already marked by
$\boxtimes)$. What remains is the diagram for $\PC(7,5)$, where now
$n<2m$. If this were transposed to get the case of $\PC(7,2)$ then
the previous procedure would take the positions indicated by
$\otimes$ into $S$, namely (1,3) and (2,4). Now one must remove the
first two columns, removing the $5\times2$ block indicated by the
``bullets'' $\bullet$, and then the first two rows, which will
remove the $2\times5$ block indicated by the small circles $\circ$.
One is now in the $\PC(5,2)$ case but the remaining columns and rows
begin with the third, so (3,5) and (4,6) are taken into $S$; these
are in the positions indicated by
$\large{\tl}\hspace{-.12in}\times$\, (an approximation to an x
included in a triangle). Now the first two of the remaining columns
are removed (columns 3 and 4 of the original), which removes the
$3\times2$ block indicated by the black triangles $\btl$, and then
first two of the remaining rows are removed (rows 3 and 4 of the
original), which removes the $2\times3$ block indicated by the
triangles $\tl$ (including those already marked
$\tl\hspace{-.12in}\times$ . What is left is a $\PC(3,1)$ where
(5,7), in the position marked by an x in a right pointing triangle
$\tgr\hspace{-.12in}\footnotesize{\times}$ is taken into $S$, and
the block marked by the black right pointing triangles $\btr$ and
the open right pointing triangles $\tgr$ are removed. This leaves
finally a $\PC(2,1)$ and the last entry taken into $S$, namely
(5,6), is marked simply $\times$.

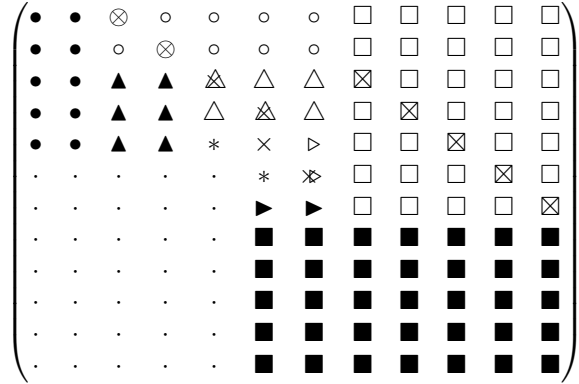
\begin{equation*}
\setcounter{MaxMatrixCols}{16}
\begin{pmatrix}
\bt&\bt&\otimes &\circ &\circ &\circ &\circ &\sq &\sq &\sq &\sq &\sq \\
\bt&\bt&\circ &\otimes &\circ &\circ &\circ &\sq &\sq &\sq &\sq &\sq \\
\bt&\bt&\btl&\btl&\tl\hspace{-.124in}\times &\tl &\tl&\boxtimes &\sq &\sq &\sq &\sq \\
\bt&\bt&\btl&\btl&\tl &\tl\hspace{-.123in}\times&\tl &\sq &\boxtimes &\sq &\sq &\sq \\
\bt&\bt&\btl&\btl&\ast &\times &\tgr &\sq &\sq &\boxtimes &\sq &\sq \\
\cdot&\cdot&\cdot&\cdot&\cdot&\ast &\tgr\hspace{-.122in}\times &\sq &\sq &\sq &\boxtimes &\sq \\
\cdot&\cdot&\cdot&\cdot&\cdot&\btr &\btr &\sq &\sq &\sq &\sq &\boxtimes \\
\cdot&\cdot&\cdot&\cdot&\cdot&\bsq&\bsq&\bsq&\bsq&\bsq&\bsq&\bsq \\
\cdot&\cdot&\cdot&\cdot&\cdot&\bsq&\bsq&\bsq&\bsq&\bsq&\bsq&\bsq \\
\cdot&\cdot&\cdot&\cdot&\cdot&\bsq&\bsq&\bsq&\bsq&\bsq&\bsq&\bsq \\
\cdot&\cdot&\cdot&\cdot&\cdot&\bsq&\bsq&\bsq&\bsq&\bsq&\bsq&\bsq \\
\cdot&\cdot&\cdot&\cdot&\cdot&\bsq&\bsq&\bsq&\bsq&\bsq&\bsq&\bsq
\end{pmatrix}
\end{equation*}
\setcounter{MaxMatrixCols}{10}

The proof that the upper triangular $S$ is Frobenius follows the
same reduction procedure as for the cyclic $S$. The blocks have been
removed in the example as they would be in the reverse induction.
(That predecessors now need not be unique causes a only a slight
complication in the proof; there are natural choices.) To compute
the $r$ matrix one shows that the subspace $L$ spanned by all
$e_{ij}, i\ne j$ in the first block removed in every pair (the black
ones in the illustration) together with all $e_s, s\in S$ is
Lagrangian. Its dual $L'$ consists of the Cartan subalgebra of
$\spl(n)$ together with all $e_{ij}, (i,j) \notin S$ in the second
block removed in every pair. Although $L'$ is generally not
Lagrangian the computation is not difficult, cf.
Section~\ref{sec:Lagrange}.

\section{Seaweed algebras}\label{sec:seaweed} Dergachev and Kirillov
\cite{DergKir:Index} define a \emph{seaweed} subalgebra $\frak k$ of
a simple Lie algebra $\g$ to be one generated by a Cartan subalgebra
together with the root spaces of some subset of the simple roots,
both positive and negative. (The suggestive name comes from the
picture of such an algebra when $\g = \spl(n)$.) Equivalently,
$\frak k$ is the intersection of a positive parabolic subalgebra
(omit some of the negative roots) and of a negative parabolic; A.
Joseph \cite{Joseph:biparaI} has therefore also called them
``biparabolic". For seaweed subalgebras of $\spl(n)$, which include
all $\PC(n,m)$, Dergachev and Kirillov define a functional $F_S$ the
importance of which is that it is always regular. As mentioned
above, The Dergachev-Kirillov functional is constructed by taking
into $S$ elements on antidiagonals starting at the corners and
proceeding until one reaches the main diagonal. For $\PC(n,m)$ it
consists of all pairs $(i,j) \in \Pi(n,m)$ with $i \ne j$ and either
$i+j = m+1, i>j$ or $i+j = n+1, i<j$ or $i+j= n+m+1, 1\ge n-m$. (The
first of these three sets will be empty for $m=1$ and the last will
be empty for $m= n-1$; the diagram for $\PC(7,3)$ is in
Section~\ref{sec:gamma}.) Since this functional is always regular,
when $(n,m) = 1$ it is Frobenius and small in our sense. While there
is yet no general closed formula for the index of a seaweed algebra,
Dergachev and Kirillov give an algorithm, using a `meander' built
from their $S$, which once the omitted roots are specified rapidly
computes the index. This meander is the graph $\gamma(S)$, which for
their $S$ is always a disjoint union of loops and chains (no
branching); the index is 2(\#loops)$+$(\#chains)$+$(\#isolated
vertices)$-$1. For $\PC(n,m)$ with $(n,m) = 1$ it is easy to verify
that the graph $\gamma(S)$ is a single chain (with possible
reversals of arrows), giving an index of zero, but $\G(S)$ generally
has branches. The automorphism of $\spl(n)$ carrying $\PC(n,m)$ to
$\PC(n,n-m)$ carries the Dergachev-Kirillov functional of the former
to that of the latter.

\section{The Belavin--Drinfel'd solutions to the MCYBE}\label{sec:MCYBE}

Belavin--Drinfel'd have given an explicit construction of all
solutions to the MCYBE associated to a simple Lie algebra $\g$ in
\cite{BelDrin:2}. We do not need a full description of their work in
what follows; the reader is referred to \cite{BelDrin:2} for
details. We shall use the fact, however, that the set of solutions
is a finite disjoint union of components each of which is determined
by an ``admissible triple'', which is a bijection between two
subsets of positive simple roots of $\g$, satisfying certain
properties.

For the case $\g= \spl(n)$ the simple roots may be identified with
the set $\{1,2,\dots, n-1\}$ and an admissible triple $\mathcal{T}$
is in effect a bijection $T:S_1\to S_2$ between subsets of
$\{1,\dots,n-1\}$ such that (1) for every $i\in S_1$ there is an $r$
with $T^r(i) \notin S_1$ and (2) $T$ preserves adjacency, i.e., if
$i,j\in S_1$ with $|i-j| = 1$ then $|T(i)-T(j)| = 1$. There is a
natural partial order among triples where $\mathcal{T}=(S_1,S_2,
T)\prec \mathcal{T}'=(S_1',S_2',T')$ if $T$ is the restriction of
$T'$ to some subset $S_1$ of $S'_1$. Any admissible triple with
$\#S_1 = n-2$ is then maximal in the partial order of triples, and
these were all determined in \cite{GG:Boundary}: Denoting the
omitted element of $S_1$ by $n-m$, it must be the case that $m$ and
$n$ are relatively prime (suggesting a relation with the fact that
$\PC(n,m)$ is then Frobenius); the omitted element of $S_2$ is $m$,
and $T$ sends every $i\in S_1$ to $i+m$ understood modulo $n$. With
the standard triangular decomposition, $\spl(n)=\n^-\oplus \h\oplus
\n^+$, each solution to the MCYBE is of the form $\gamma + \beta
+\alpha$ where $\beta \in \h\wedge\h$ and $\alpha\in \n^+\wedge
\n^-$ are determined by $\mathcal{T}$ and $\gamma =
\sum_{i<j}e_{ij}\wedge e_{ji}$. In particular $\alpha$ is uniquely
determined by $\mathcal{T}$ and so we write it as
$\alpha(\mathcal{T})$.

Let $r'=\gamma +\beta +\alpha(\mathcal{T'})$ be a solution to the
MCYBE associated to some triple $\mathcal{T'}$. If $\mathcal{T}\prec
\mathcal{T'}$, then it was asserted in \cite{GGS:Construction} that
$r=\gamma +\beta+\alpha(\mathcal{T})$ is in the closure of the orbit
of $r'$ under the operation of $SL(n)$. We show this here. (The idea
of the proof probably works for all simple $\g$.) If $h \in \h$ then
$[h, \gamma] = [h,\beta] = 0$, so it is sufficient to find an $h$
such that $\exp(th)\alpha(\mathcal{T'})\exp(-th) \to
\alpha(\mathcal{T})$ as $t \to \infty$. Now $\mathcal{T}'$
establishes a partial order on the roots which we can extend to a
linear order, and by renumbering we may suppose that this coincides
with the natural order $1,\dots,n-1$. If $h\in\h$ has diagonal
entries $h_1,\dots,h_n$ then setting $\lambda_i = h_i-h_{i+1}, i =
1,\dots,n-1$ one has $[h, e_{i,i+1}]=\lambda_ie_{i,i+1}$. Since
$\alpha'$ will now be a sum of elements of the form $e_{i,i+1}
\wedge e_{j,j+1}$ with $j > i$,
$\exp(th)\alpha(\mathcal{T'})\exp(-th)$ will have a finite limit as
$t \to \infty$ whenever $\lambda_1 \le \lambda_2 \le \dots \le
\lambda_{n-1}$. To make the limit agree with $\alpha(\mathcal{T})$
it is sufficient to make some of these inequalities strict, and it
is trivial to find an $h\in \h$ with these properties.
Interestingly, in all examples studied, if we begin with $P(n,m)$,
let $F$ be the cyclic functional and take $h =\hat F$, the
associated principal element, then this $h$ has the required
properties with respect to the associated $\mathcal{T'}$. Moreover,
for the principal element $h$, we have
$\exp(th)\alpha(\mathcal{T'})\exp(-th) \to \alpha(\mathcal{T})$ as
$t\to\infty$ where $T$ omits from $T'$ only certain mappings sending
an $i$th root to a $j$th root with $j < i$. In general it does not
remove all of these except when there is just one. For example, with
$P(5,2)$ where the progression of roots is $1\to 3\to 2 \to 4$ the
map from $3$ to $2$ is removed, this being the only mapping in
``reverse direction". However, with $P(8,5)$ where the progression
of roots is $5 \to 2 \to 7 \to 4 \to 1 \to 6 \to 3$, only the
mappings $5 \to 2$ and $4 \to 1$ are removed while $7 \to 4$ and $6
\to 3$ remain.

\section{Local rings associated to a graph; a reconstruction
theorem}\label{sec:local} This section is added to show that
questions about graphs can be expressed as ones in local algebra
(but that does not necessarily make them any easier). To every graph
we associate two local rings, ``full'' and ``reduced''. The first
captures all the information in the graph with but one exceptional
case.

Suppose that we have a graph $\G$ whose sets of edges and vertices
will be denoted by $E$ and $V$, respectively. In the polynomial ring
$K[E]$ generated by the edges, let $I$ be the ideal generated by all
relations of the form $e_1e_2 = 0$ whenever $e_1, e_2 \in E$ have a
common vertex and set $K\G = K[E]/I$; we will call it the \emph{full
local ring of $\G$}. Since $I$ is homogeneous, $K\G$ continues to be
graded. Denote its radical (augmentation ideal) by $J$. Then $J$ is
nilpotent, and its index of nilpotence is $\mn(\G) +1$. The
dimension of $J^k$ is the number of disjoint $k$-tuples of mutually
disjoint edges in $\G$. For example, if $\G$ is a square then $\dim
J = 4$ (the number of edges) and $\dim J^2 = 2$. One also has
$\mn(\G) = 2$ and $J^3 = 0$. Knowing $K\G$ does not determine $\G$
since the triangle and three-pointed star both have local ring
isomorphic to a ring $\mathcal{R}_3$ generated over $K$ by three
variables all of whose squares and products vanish, but we will show
that amongst connected graphs this is the only counterexample. The
triangle and the three pointed star are clearly the only graphs with
$\mathcal{R}_3$ as local ring since the radical has dimension three,
so there are exactly three edges, and the three-pointed star and
triangle are the only configurations in which each has a vertex in
common with every other. The ring can not recognize that the
triangle has only three vertices while the star has four.  One can
tell from $K\G$ if $\G$ is connected. Call a local ring \emph{graph
connected} if its radical $J$ can not be written as a direct sum of
subspaces $J_1$ and $J_2$ such that $x_1 \in J_1, x_2 \in J_2$ and
$x_1, x_2 \ne 0$ imply $x_1x_2 \ne 0$. (There are other concepts of
connectedness.) It is evident that a graph is connected if and only
if its full local ring is graph connected.

\begin{Theorem} A connected graph $\G$ can be reconstructed from
$K\G$ except when $K\G \cong \mathcal{R}_3$; connected graphs with
isomorphic full local rings are isomorphic except for the triangle
and the three-pointed star.
\end{Theorem}

\noindent\textsc{Proof.} We show that one can reconstruct a graph
$\G$ from its full local ring $K\G$ as long as the ring is graph
connected and the dimension of its radical $J$ is not three. The
cases where $\dim J$ are smaller than three are trivial, so we may
assume that it is at least four. A \emph{zero algebra} is one in
which all products vanish. If the graph $\G$ contained no triangles
then we could reconstruct it easily from $K\G$. The vertices would
then be in one-one correspondence with the maximal zero subalgebras
of $K\G$ and the vertices represented by two zero subalgebras $Z$
and $Z'$ would then be joined by an edge exactly when $Z\cap Z' \ne
0$, in which case the intersection would have dimension one. (Note
that having $Z\cap Z' \ne 0$ is equivalent to the existence of a
non-zero element which annihilates both, in which case there must
also be non-zero elements in $Z$ and $Z'$, respectively, whose
product is not zero, else $Z$ would coincide with $Z'$.) The only
problematic case is that where the maximal zero algebra $Z$ has
dimension three, i.e., is isomorphic to the radical of
$\mathcal{R}_3$, for then we can not tell immediately if it has come
from a triangle which is a subgraph of $\G$ or a vertex whose star
is the three pointed star. Since the dimension of the radical is
greater than three there must be some maximal zero subalgebra $Z'$
with $Z\cap Z' \ne 0$, else the local algebra would not be graph
connected. For in the notation above we could then take $J_1 = Z$
and $J_2$ to be the sum  of all other zero subalgebras of the
radical. If $Z\cap Z'$ has dimension one then $Z$ and $Z'$ represent
vertices. The only other possibility is that the dimension is two,
in which case $Z$ came from a triangle and $Z'$ from one of its
vertices; that vertex is determined by the particular two
dimensional subspace $Z\cap Z'$ of $Z$. $\Box$
\bigskip

Let $\bigwedge V$ be the exterior algebra generated by the vertices
$V$ of $\G$ and $\bigwedge_{\mathrm{even}}$ be its even subalgebra,
i.e., the subalgebra generated by the unit element and all elements
of the form $v\wedge v'$. Now assign an arbitrary orientation to
every edge $e$ of $\G$ so that we can say which vertex is initial
and which terminal. Then we can define a ring morphism $K\G \to
\bigwedge_{\mathrm{even}}$ by sending every edge $e$ with initial
vertex $v$ and terminal vertex $v'$ to $v\wedge v'$. The image
$(K\G)_{\mathrm{red}}$ is the \emph{reduced local ring}. It does not
depend on the choice of orientation, is naturally graded, and is
again a commutative local ring whose radical still has index of
nilpotence equal to $\mn(\G) +1$, but its dimension over $K$ may be
smaller than that of $K\G$. For example, if $\G$ is a square with
sides (in counterclockwise order) $x,y,z,w$ and respective vertices
$a,b,c,d$, then except in the case of characteristic 2, $J^2$ has
dimension 2, being spanned by $xz$ and $yw$, while its image has
dimension only 1, being spanned by $a\wedge b\wedge c\wedge d$.

\end{document}